\documentclass[11pt]{article}

\usepackage{amsmath, amsthm, amssymb}
\usepackage{hanging}

\newtheorem{definition}{Definition}

\newcommand{\Sys}{\mathcal{S}}
\newcommand{\SRE}{\mathcal{S}_{\text{\tiny RE}}}
\newcommand{\SLBI}{\mathcal{S}_{\text{\tiny LBI}}}
\newcommand{\TRE}{T_{\text{\tiny RE}}}

\newcommand{\LBI}{(x \lor \lnot x) \to y}

\title{A Note on the Incompleteness of\\
        Gödel's Incompleteness Theorems}
\author{Jeffrey Uhlmann \\
Dept.\ of Electrical Engineering and Computer Science\\
University of Missouri - Columbia}
\date{}

\begin{document}

\maketitle

\begin{abstract}
In this note we observe that automated theorem provers
(ATPs) that recursively enumerate theorems in a particular
way will fail to identify some valid theorems for a reason 
that is analogous to how Gödel proved the existence of 
what are now referred to as Gödel statements. This 
observation has no significant practical or theoretical
implications, but it may be of pedagogical value for
honing the intuition of students about recursive 
enumeration in the context of ATP.
\end{abstract}
\noindent Keywords: Gödel's Incompleteness Theorems,
Recursive Enumeration, Law of Excluded Middle, Automated
Theorem Proving.

\section{Introduction}

In this note we observe that automated theorem provers
(ATPs) that recursively enumerate theorems in a particular
way will fail to identify some valid theorems for a reason 
that is analogous to how Gödel proved the existence of 
what are now referred to as Gödel statements \cite{goedel1}. 
This observation has no significant practical or theoretical
implications, but it may be of pedagogical value for
honing the intuition of students about recursive 
enumeration in the context of ATP.

The structure of this paper is as follows. We begin 
by informally discussing an inference rule based on the 
Law of Excluded Middle (LEM). We then show that this
rule, LEM-Based Inference (LBI), represents an obstacle
to the bottom-up enumeration of theorems by the
application of a given set of inference rules to a given
set of axioms.  

\section{LEM-Based Inference}

Use of the Law of Excluded Middle (LEM) as a basis for
inference has the following form:
\begin{equation}
   (x \lor \lnot x) \to y,
\end{equation}
where the presumption is that either $x$ or $\lnot x$ must be
true, therefore either $(x\to y)$ or $(\lnot x\to y)$ must be true. 
This bivalence was stated explicitly by Aristotle in the 
{\em Organon} \cite{organon}, and its use as an inference rule has a 
long history. Notable examples include the following:
\begin{itemize}
\item In 1914, Littlewood proved a result relating to the
Prime Number Theorem based on its implication from the truth
of the Riemann Hypothesis and from its negation \cite{littlewood1,littlewood2}. 
\item In 1918 up to the 1930s, a collection of theorems relating to the 
Gauss Class Number Conjecture were proven via LBI relative 
to the Generalized Riemann Hypothesis \cite{ir1990} 
\item A series of results relating to Stone-\v{C}ech 
compactification were proven via LBI relative to the Continuum 
Hypothesis \cite{vanmill}. 
\item In 1975, Erdos and Nicolas proved a bound on Euler's totient
function using LBI relative to the Riemann Hypothesis \cite{nicolas1,nicolas2}.
\item Fermat's Last Theorem (FLT) was proven via LBI \cite{lipblog} in 1995
by Andew Wiles \cite{wiles} with contributions by Robert Taylor \cite{taylor}.
\end{itemize}
LBI-based proofs are also be found in the 
engineering and computer science literature and
are similarly not regarded as controversial.

\section{LBI, RE, and Formal Systems}

We begin with the following definitions.

\begin{definition}
{\em Axiomatic System (AS)}: An axiomatic system $\Sys$ consists of a set of axioms, along with a set of inference rules from which a set of provable theorems can be derived.
\end{definition}

\begin{definition}
{\em Recursively Enumerable (RE) Axiomatic System}: An axiomatic formal system $\SRE$ is said to be RE if its set of theorems $\TRE$ is recursively enumerable, i.e., there exists a Turing machine that can enumerate all $t \in \TRE$.
\end{definition}

\begin{definition}
{\em LBI-Accepting System}: An axiomatic system $\SLBI$ is said to be LBI-accepting if for any formula $y$ and any proposition $x$, if $\LBI$ is a valid inference, then $y$ is a theorem in the system.
\end{definition}

The complete set of theorems for a given system can be recursively enumerated by examining every proposition (i.e., possible sequence of symbols) and assessing whether it can be proven from the system's set of axioms and rules of inference. An intuitively attractive approach for performing automated theorem proving would be to perform a bottom-up ordered construction of a tree of theorem derivations starting from the system's set of axioms. Thus, for every theorem there would exist a branch representing the application of a sequence of inference rules to a sequence of theorems established earlier in the construction that gives a proof of the theorem. This tree construction would not be unique because it depends on arbitrary choices such as the specific order in which axioms, rules, and theorems are processed, but it would be expected to recursively enumerate all theorems provable within the system. We show that this is not the case.\\
~\\
\noindent {\bf Main Theorem}:
{\em The theorem set for an LBI-accepting system $\SLBI$ cannot generally be recursively enumerated using bottom-up enumeration}.\\
~\\
\noindent {\bf Proof}:
Assume that the theorem set of $\SLBI$ is RE. We know from Gödel that if $\SLBI$ is sufficiently powerful there will exist undecidable $x$ such that neither $x$ nor $\lnot x$ are elements of the RE theorem set, i.e., the truth value of $x$ cannot be established by the system. Let $y$ be an LBI-accepted theorem via $(x \lor \lnot x) \to y$, then $y$ cannot be an element of the bottom-up enumeration because neither $x$ nor $\lnot x$ exist as theorems from which $y$ can be derived.\\ 
{\small\em QED}.\vspace{11pt}

\section{Conclusion}

We have given a brief explanation for how the basic approach underpinning Gödel's incompleteness theorems can be applied in an intuitively-relatable context to show the ``incompleteness'' of a bottom-up method for enumerating the set of theorems provable by a given set of axioms and rules of inference.

\bibliographystyle{amsplain}

\begin{thebibliography}{99}

\bibitem{goedel1}
Kurt Gödel, ``Über formal unentscheidbare Sätze der Principia Mathematica und verwandter Systeme, I,'' {\em Monatshefte für Mathematik und Physik}, v.\ 38 n.\ 1, pp.\ 173–198, 1931.

\bibitem{organon} 
Aristotle, {\em Organon}, (Andronicus of Rhodes, ed.), 9 October (CDT), 40\,BC.


\bibitem{littlewood1}
Littlewood, J.E.,``Sur la distribution des nombres premiers,'' {\em Comptes Rendus}, 158, 1914.

\bibitem{littlewood2}
Ingham, A.E, {\em The Distribution of Prime Numbers, Cambridge Tracts in Mathematics and Mathematical Physics}, vol.\ 30, Cambridge University Press, 1932. 

\bibitem{nicolas1}
Erdős, P.\ and Nicolas, J.L., ``Répartition des nombres superabondants, {\em Bull.\ Soc.\ Math.\ France}, 79 (103): 65–90, 1975.

\bibitem{nicolas2}
Sárközy, A., ``Jean-Louis Nicolas and the partitions,'' {\em The Ramanujan Journal}, 9 (1–2): 7–17, 2005.

\bibitem{ir1990}
Ireland, Kenneth and Rosen, Michael, {\em A Classical Introduction to Modern Number Theory} (Second edition), New York: Springer, 1990.

\bibitem{vanmill}
van Mill, Jan, ``An introduction to $\beta\omega$,'' in Kunen, Kenneth; Vaughan, Jerry E.\ (eds.), {\em Handbook of Set-Theoretic Topolog}y, North-Holland, pp.\ 503–560, 1984.

\bibitem{lipblog}
Lipton, RJ, ``Proof Checking: Not Line by Line,'' {\em Gödel's Lost Letters and P=NP} (blog), 13 June, 2020.

\bibitem{wiles}
Wiles, Andrew, ``Modular elliptic curves and Fermat's Last Theorem,'' {\em Annals of Mathematics}, 141 (3): 443–551, 1995.

\bibitem{taylor}
Taylor R, Wiles A, ``Ring theoretic properties of certain Hecke algebras,'' {\em Annals of Mathematics}, 141 (3): 553–572, 1995.




\end{thebibliography}

\end{document}